\input amstex
\documentstyle{amsppt}
\NoBlackBoxes
\nologo \magnification1200 \pagewidth{6 true in}
\pageheight{9 true in} 
\topmatter
\title Extreme values of $\arg L(1,\chi)$
\endtitle
\author{ Youness Lamzouri }
\endauthor
\address School of Mathematics,
Institute for Advanced Study,
1 Einstein Drive,
Princeton, NJ, 08540
USA
\endaddress
\email{Lamzouri{\@}math.ias.edu}
\endemail
\thanks
\endthanks
\endtopmatter
\define \ex{\Bbb E}
\define \La{\Cal{L}_q}
\document
\head {1. Introduction }
\endhead
\noindent
Values of Dirichlet $L$-functions at $s=1$ have always attracted great attention due to their central role in number theory. In particular the prime number theorem for arithmetic progressions relies on the non-vanishing of $L(1,\chi)$ for any non-principal character $\chi$  $(\text{mod  } q)$. Moreover improving the existing lower bounds for $|L(1,\chi)|$ (the famous Siegel's bound $|L(1,\chi)|\gg_{\epsilon}q^{\epsilon}$) would imply many important consequences.

Various authors have studied the distribution of these values. (one can refer to [5] for a history of the subject including results for other families of $L$-functions). In [3] A. Granville and K. Soundararajan studied the distribution of extreme values of $|L(1,\chi)|$, showing that the proportion of characters $\chi \ (\text{mod }q)$ for which $|L(1,\chi)|>e^{\gamma}\tau$ equals
$$ \exp\left(-\frac{e^{\tau-C_1-1}}{\tau}(1+o(1))\right),$$
uniformly in the range $1\ll \tau\leq \log\log q-20$, where
$$C_1:=\int_0^1\log I_0(t)\frac{dt}{t^2} + \int_1^{\infty} (\log
 I_0(t)-t)\frac{dt}{t^2},\tag{1}$$
 and $I_0(t): = \sum_{n=0}^{\infty}
(t/2)^{2n}/n!^2$ is the modified Bessel function of order $0$. In [6] the author studied the distribution of the values $L(1,\chi)$ in the complex plane by considering the joint distribution function of $|L(1,\chi)|$ and $\arg L(1,\chi)$ when the norm is large and the argument is bounded, where $\arg L(1,\chi)$ is defined by continuous variation
along the real axis from infinity taking the argument at infinity to be zero.

In this paper we are concerned with the study of the distribution of extreme values of $\arg L(1,\chi)$, as $\chi$ varies over primitive characters modulo a large prime $q$. Refining ideas of J.E. Littlewood [7] we first show that assuming the Generalized Riemann Hypothesis we have that (see Corollary 2.6 below)
$$|\arg L(1,\chi)|\leq \log\log\log q +C_2 +\log 2 +o(1),\tag{2}$$
where
$$ C_2:= \lim_{x\to\infty}\left(\sum_{p\leq x}\arctan\left(\frac{1}{\sqrt{p^2-1}}\right)-\log\log x\right)=0.2937504<\gamma.\tag{3}$$
Moreover we can exhibit extreme values of $\arg L(1,\chi)$ from Theorem 1 below
\proclaim{Corollary 1} For any $\epsilon>0$ there exists $q(\epsilon)>0$ such that if $q$ is a prime number  and $q\geq q(\epsilon)$, then  for $\delta\in \{-1,1\}$, there are at least $q^{1-1/(\log q)^{\epsilon}}$ non-principal characters $\chi\mod q$ for which
$$ \delta\arg L(1,\chi)\geq \log\log\log q+ C_2-\epsilon.\tag{4}$$
\endproclaim
Therefore we can see that the difference between this result and the conditional bound (2) is the constant $\log 2$, and we can ask ourselves which of these two results correspond to the true nature of extreme values of $\arg L(1,\chi)$. This gives another motivation to study the behavior of the distribution function
 $$\Psi_q(\tau):=\frac{1}{\phi(q)}|\{\chi \ (\text{mod } q), \chi\neq\chi_0 : \arg L(1,\chi)>\tau\}|, \text{ for } \tau>0.$$

{\it Remark 1.} Concerning the function 
$$\Phi_q(\tau):=\frac{1}{\phi(q)}|\{\chi \ (\text{mod } q), \chi\neq\chi_0 : \arg L(1,\chi)<-\tau\}|,$$
for $\tau>0$, one can observe that 
$$\Phi_q(\tau)=\Psi_q(\tau),$$
simply because the complex characters $\chi \ (\text{mod }q)$ occur in pairs, and for these we have $\arg L(1,\overline{\chi})=-\arg L(1,\chi)$  (for real characters $\chi$ we have $\arg L(1,\chi)=0$ by the Euler product representation of $L(1,\chi)$). Therefore all the results obtained for $\Psi_q(\tau)$ hold for $\Phi_q(\tau)$.

In [2] P.D.T.A Elliott studied the distribution of $\arg L(s,\chi)$ in the half-plane Re$(s)>1/2$ . Although he considered a larger family (all the non-principal characters to prime moduli not exceeding $Q$), his methods shows that $$\lim\Sb q\to\infty\\ q \text{
 prime }\endSb\Psi_q(\tau)=\Psi(\tau)$$ exists and is a continuous function of $\tau$. Moreover he indicated how to compute its characteristic function. In his Ph.D thesis, W.R. Monach [8] showed that $1-\Psi(\tau)$ is the distribution function of a sum of independent random variables, and used this fact to estimate the tail of this distribution. Indeed he proved that for $\tau\geq \frac{\pi}{2}$, there exist positive constants $a_1,a_2>0$ for which
 $$ \exp\left(-e^{a_1e^{\tau}}\right)\leq \Psi(\tau)\leq \exp\left(-e^{a_2e^{\tau}}\right).$$
This indicates that $\Psi_q(\tau)$ should decay ``triple exponentially'' as $\tau\to\infty$. Our first Theorem estimates the function $\Psi_q(\tau)$ in a wide range of $\tau$ (which we expect is  the full range) and confirms this conjecture. Moreover it also improves on Monach's bounds for $\Psi(\tau)$. 
\proclaim{ Theorem 1}
Let $q$ be a large prime. Uniformly for
 $1\ll \tau\leq \log\log\log q+C_2-o(1)$ we have
$$ \Psi_q(\tau)=\exp\left(-\frac{\exp\left(e^{\tau-C_2}-C_1-1\right)}
{e^{\tau-C_2}}
\left(1+O\left(\frac{1}{e^{\tau/2}}\right)\right)\right).$$
The same asymptotic also holds for $\Psi(\tau)$ but for arbitrary $\tau$.
\endproclaim
{\it Remark 2.} If the asymptotic for $\Psi_q(\tau)$ were to persist to a slightly bigger range $1\ll \tau\leq \log\log\log q+C_2+\epsilon$ for any $\epsilon>0$ then this implies the conjecture that
$$ \align
&\max\Sb \chi \ (\text{mod } q)\\ \chi\neq\chi_0\endSb \arg L(1,\chi)=\log\log\log q+ C_2+o(1), \text{ and }\\
& \min\Sb \chi \ (\text{mod } q)\\ \chi\neq\chi_0\endSb \arg L(1,\chi)=-\log\log\log q- C_2+o(1).\\
\endalign
$$

To establish Theorem 1, the main idea consists of relating the Laplace transform of $\arg L(1,\chi)$ with large purely imaginary moments of $L(1,\chi)$ (in this case this approach works better than estimating the moments of $\arg L(1,\chi)$). Then we evaluate these imaginary moments (in our case the moments we need are $M_q(-is,is)$ where $s$ is large positive number, see (5) below for the definition), and use the saddle point method to recover the asymptotic for the distribution function. Therefore the next step is to find an asymptotic formula for complex moments of $L(1,\chi)$.

Let us first define some notation which will be used throughout this paper. For $z$ a complex number, we define the ``$z$th divisor function'' $d_z(n)$, to be the multiplicative function such that $d_z(p^a)=\Gamma(z+a)/\Gamma(z)a!$, for any prime $p$ and any integer $a\geq 0$. Then $d_z(n)$ is the coefficient of the Dirichlet series $\zeta(s)^z$ for Re$(s)>1$. Furthermore
$\log_j$ will denote  the $j$-th iterated logarithm, so that $\log_1n
=\log n$ and $\log_j n=\log (\log_{j-1}n) $ for each $j\geq 2$.
Let $$S_q:= \{\chi  \ (\text{mod }q), \chi\neq \chi_0, \chi \text{ is non-exceptional} \},$$
where we define a character $\chi$ to be {\it exceptional } if there exists $s$ with $\text{Re}(s)\geq 1-c/(\log (q(\text{Im}(s)+2)))$ and $L(s,\chi)=0$, for some sufficiently small constant $c>0$. One expects that there is no such a character, but what is known
unconditionally (see [1]), is that these characters if they exit, must be very rare.
Indeed each $\chi$ must be real, and between any
two powers of $2$ there is at most one fundamental discriminant $D$
such that $\chi=\left(\frac{D}{\cdot}\right)$. Since $q$ is prime in our case, there is at most
one such exceptional character of conductor $q$.

Granville and Soundararajan (unpublished) proved an asymptotic formula for the moments
$$ M_q(z_1,z_2):= \frac{1}{\phi(q)}\sum_{\chi\in S_q}L(1,\chi)^{z_1}L(1,\overline{\chi})^{z_2},\tag{5}$$
where $z_1,z_2$ are complex numbers satisfying $|z_1|,|z_2|\leq \log q/(\log_2q)^3.$  Using their method, one can improve this range to $\log q/(50(\log_2q)^2)$ (see Theorem 9.2 of [6]). Their idea consists of using an induction on certain sums of divisor functions to control the off-diagonal terms of the moments. Using a different approach (based on zero density estimates for $L(s,\chi)$) we improve this range to
$$ R(q):=\frac{\log q\log_4q}{10\log_2 q\log_3 q}.$$
\proclaim{Theorem 2}
Let $q$ be a large prime. Then uniformly
for all complex numbers $z_1,z_2$ in the region
 $|z_1|,|z_2|\leq R(q)$,
we have
$$ M_q(z_1,z_2) =
\sum_{n=1}^{\infty}\frac{d_{z_1}(n)d_{z_2}(n)}{n^2}+O\left(\exp\left(-\frac{\log q\log_4q}{50\log_2 q}\right)\right).\tag{6}$$
\endproclaim
 As in [5] and [6] we can compare the distribution of $\arg L(1,\chi)$ to that of an appropriate probabilistic model. Let $\{X(p)\}_{p\text{ prime }}$ be independent random variables uniformly distributed on the unit circle, and define the random Euler products $L(1,X):=\prod_p(1-X(p)/p)^{-1}$ (these products converge with probability $1$). Indeed one can observe that the main term on the RHS of (6) corresponds to  ${\Bbb
E}\left(L(1,X)^{z_1}\overline{L(1,X)}^{z_2}\right)$ (this follows from Lemma 3.2 below), and that $\Psi(\tau)=\text{Prob}(\arg L(1,X)>\tau)$ (this has been proved in [8]).

\head{2. Bounds for $\arg L(1,\chi)$ and $\log L(s,\chi)$}\endhead
\noindent We begin by proving the following useful Lemma which will be used later in several places.
\proclaim{Lemma 2.1} Let $x$ be a large number, and  $0\leq \alpha\leq\log_3x/(4\log x)$. Then for $1\leq \sigma<4$ we have
$$\sum_{p\leq x^{\sigma}}\frac{1}{p^{1-\alpha}}=\log\log x (1+o(1)).$$
\endproclaim
\demo{Proof}The lower bound follows from the prime number theorem which implies that
$$ \sum_{p\leq x^{\sigma}}\frac1p =\log\log x + O(1).$$
To prove the upper bound we split the sum into two parts. In the range $p\leq e^{1/\alpha}$ we use the estimate $p^{\alpha}=1+O(\alpha\log p)$ and the asymptotic formula $\sum_{p\leq y}\log p/p=\log y+O(1)$, while in the range $e^{1/\alpha}\leq p\leq x^{\sigma}$ we use the trivial bound $p^{\alpha}\leq (\log_2x)^{\sigma/4}$. Thus we get
$$
\align
\sum_{p\leq x^{\sigma}}\frac{1}{p^{1-\alpha}}&\leq \sum_{p\leq e^{1/\alpha}}\frac{1+O(\alpha\log p)}{p} +(\log_2x)^{\sigma/4}\left(\log\left(\frac{\log x^{\sigma}}{\log e^{1/\alpha}}\right)+O(1)\right)\\
&\leq \log_2 x-\log_4 x + (\log_2x)^{\sigma/4} \left(\log_4 x +O(1)\right)=\log\log x(1+o(1)),\\
\endalign
$$
which completes the proof.
\enddemo
In the next Lemma we establish the classical bound for $\log L(s,\chi)$ on the line Re$(s)=1$ if $\chi\neq\chi_0$ is a non-exceptional character.
\proclaim{Lemma 2.2} If $\chi\neq\chi_0$ is a non-exceptional character $\mod q$, then for all $t\in \Bbb{R}$ we have
$$ |\log L(1+it,\chi)|\ll \log_2 q(|t|+2).\tag{2.1}$$
\endproclaim
\demo{Proof}
Consider the circles with center $s_0=1+1/(\log q(|t|+2))+it$ and
radii $r:=1/(\log q(|t|+2))<R:=(1+c)/(\log q(|t|+2))$, for an appropriately small constant $c>0$. Then the smaller circle passes
through $1+it$. From our assumption on $\chi$ we may choose $c$ such that $\log L(z,\chi)$ is analytic inside the
larger circle, by the classical zero free region of $L(z,\chi)$. For a point $z$ on the larger circle we have the classical
bound $\log|L(z,\chi)|\leq \log_2q(|t|+2)+O(1)$. Thus by the
Borel-Caratheodory Theorem we deduce that
$$\align
|\log L(1+it,\chi)|&\leq \frac{2r}{R-r}\max_{|z-s_0|=R} \text{Re} (\log
L(z,\chi)) +\frac{R+r}{R-r}|\log
L(s_0,\chi)|\\
&\ll\log_2q(|t|+2),\\
\endalign
$$
using that $|\log
L(s_0,\chi)|\leq \log\zeta(1+1/(\log q(|t|+2)))=\log_2q(|t|+2)+O(1).$
\enddemo
From this Lemma we can deduce the classical bound for $\arg L(1,\chi)$
\proclaim{Corollary 2.3} If $\chi\neq\chi_0\mod q$, then
$$\arg L(1,\chi)\ll \log_2 q.$$
\endproclaim

The main ingredient to establish these classical bounds is the zero free region for $L(s,\chi)$. Therefore any improvement will depend on our knowledge of the location of its zeros. Indeed we will improve the bound (2.1) (see Corollary 2.5 below) for characters $\chi \ (\text{mod } q)$ for which $L(s,\chi)$ has no zeros in a larger region inside the critical strip (this assumption is true for almost all characters by the classical zero density result (2.2)). In particular we obtain the bound (2) for $\arg L(1,\chi)$ under the GRH. The main ingredient to establish these results is the following Lemma which corresponds to Lemma 8.2 of [4]
\proclaim{Lemma 2.4} Let $s=\sigma+it$ with $\sigma>1/2$ and $|t|\leq 2q$. Let $y\geq 2$ be a real number, and let $1/2\leq \sigma_0<\sigma$. Suppose that the rectangle $\{z:\sigma_0<\text{Re}(z)\leq 1, |\text{Im}(z)-t|\leq y+3\}$ contains no zeros of $L(z,\chi)$. Put  $\sigma_1=\min(\frac{\sigma+\sigma_0}{2},\sigma_0+\frac{1}{\log y})$. Then
$$ \log L(s,\chi)=\sum_{n=2}^y\frac{\Lambda(n)\chi(n)}{n^s\log n}+O\left(\frac{\log q}{(\sigma_1-\sigma_0)^2}y^{\sigma_1-\sigma}\right).$$
\endproclaim
\proclaim{Corollary 2.5} Let $\eta=\log_4q/(4\log_2q)$. Assume that $L(z,\chi)$ has no zeros in the rectangle $\{z: 5/8\leq \text{Re}(z)\leq 1\text{ and } |\text{Im}(z)|\leq 2\log^3q\}.$ Then for any $s=\sigma+it$ with $1-\eta\leq \sigma\leq 1$ and $|t|\leq \log^3q$ we have
$$ |\log L(s,\chi)|\leq \log_3 q(1+o(1)).$$
Furthermore this holds for all but at most $q^{10/11}$ characters $\chi \ (\text{mod q})$.
\endproclaim
\demo{Proof}
We use Lemma 2.4 with $1-\eta\leq \sigma\leq 1$, $\sigma_0=5/8$, $y=\log^3q$ and $\sigma_1=5/8+1/\log y$. Therefore if $L(z,\chi)$ has no zeros in the rectangle $\{z: 5/8\leq \text{Re}(z)\leq 1\text{ and } |\text{Im}(z)|\leq 2\log^3q\},$ then

$$
\align
|\log L(s,\chi)|&=\left|\sum_{n=2}^{\log^3q}\frac{\Lambda(n)\chi(n)}{n^s\log n}\right|+O\left(\frac{1}{\log^{1/9}q}\right)=\left|\sum_{p\leq\log^3q}\frac{\chi(p)}{p^s}\right|+O(1)\\
&\leq \sum_{p\leq\log^3q}\frac{1}{p^{1-\eta}}+O(1)\ll \log_3 q(1+o(1)),
\endalign
$$
by Lemma 2.1, taking $x=\log q$ there. Finally the last statement follows by taking $\sigma=5/8$ and $T=2\log^3q$ in the following zero density result of H.L. Montgomery [9] which states that for $q\geq 1$, $T\geq 2$ and $1/2\leq \sigma \leq 4/5$ we have
$$ \sum_{\chi \ (\text{mod } q)} N(\sigma,T, \chi)\ll (qT)^{3(1-\sigma)/(2-\sigma)}(\log qT)^9,\tag{2.2}$$
where $N(\sigma,T, \chi)$ denotes the number of zeros of $L(s,\chi)$ such that Re$(s)\geq \sigma$ and $|\text{Im}(s)|\leq T$.
\enddemo

Now we prove the bound (2) for $\arg L(1,\chi)$ under the Generalized Riemann Hypothesis
\proclaim{Corollary 2.6} Assume GRH. Then
$$|\arg L(1,\chi)|\leq \log_3 q +C_2 +\log 2 +o(1).$$

\endproclaim

\demo{Proof} By Lemma 2.4, the GRH implies that for $y=(\log^2 q)(\log_2 q)^6$ we have
$$ \log L(1,\chi)=\sum_{n=2}^y\frac{\Lambda(n)\chi(n)}{n\log n} +O\left(\frac{1}{\log\log q}\right).$$
Therefore extracting the imaginary parts from both sides we get
$$ \arg L(1,\chi)=\sum_{p\leq y}\sum_{n=1}^{\infty}\frac{\text{Im}(\chi(p^n))}{p^nn}+ O\left(\frac{1}{\log\log q}\right).$$
Now $\text{Im}(\chi(p^n))=\sin(n\arg\chi(p)).$ Then by the proof of Lemma 3.5 below we know that for all $\theta\in [-\pi,\pi]$
$$ \left|\sum_{n=1}^{\infty}\frac{\sin(n\theta)}{p^nn}\right|\leq \arctan\left(\frac{1}{\sqrt{p^2-1}}\right).$$
Thus we deduce that
$$ |\arg L(1,\chi)|\leq \sum_{p\leq y}\arctan\left(\frac{1}{\sqrt{p^2-1}}\right) +o(1)=\log_3 q+C_2+\log 2+o(1).$$
\enddemo

\head{3. Estimates for sums over divisor functions  }\endhead
We begin by collecting some useful estimates for the Bessel function $I_0(t)$.
\proclaim{ Lemma 3.1}
 $\log I_0(t)$ is a differentiable function with bounded derivative on $[0,+\infty)$
and  satisfies $$
  \ \ \ \log I_0(t)=\left\{\aligned &O\left(t^2\right) \ \ \ \ \text{ if } 0\leq t<1\\  &t+O\left(\log (t+1)\right) \ \text{ if } 1\leq t. \endaligned \right.
$$
\endproclaim
\demo{Proof}
The first estimate follows from the Taylor series expansion

\noindent $I_0(t)= \sum_{n=0}^{\infty}
(t/2)^{2n}/n!^2.$
For the second we use an integral representation of $I_0(t)$
$$ I_0(t)=\frac{1}{\pi}\int_{0}^{\pi}e^{t\cos\theta}d\theta\leq e^t.$$ Furthermore taking $\epsilon=\frac{1}{t}$ we deduce that
$$ I_0(t)\geq \frac{1}{\pi}\int_{0}^{\epsilon}e^{t\cos\theta}d\theta\geq \frac{\epsilon}{\pi}e^{t\cos\epsilon}\geq \frac{e^t}{10\pi t},$$
from which the second estimate follows.
Finally since $I_0(t)$ is a positive smooth function on $[0,+\infty)$ then $\log I_0(t)$ is differentiable and we have
$$ |(\log I_0(t))'|= \left|\frac{\int_{0}^{\pi}\cos\theta e^{t\cos\theta}d\theta}{\int_{0}^{\pi} e^{t\cos\theta}d\theta}\right|\leq 1.$$
\enddemo

Now we recall some easy bounds for the divisor function $d_z(n)$. First we have
$$ |d_z(n)|\leq d_{|z|}(n)\leq d_k(n),$$
for any integer $k\geq |z|.$ Furthermore
for $j\in {\Bbb N}$, and $X>3$ we have that

\noindent  $d_j(n)e^{-n/X}\leq
e^{j/X}\sum_{a_1...a_j=n}e^{-(a_1+...+a_j)/X}$, and so
$$ \sum_{n=1}^{\infty}\frac{d_j(n)}{n}e^{-n/X}\leq \left(e^{1/X}
\sum_{a=1}^{\infty}\frac{e^{-a/X}}{a}\right)^j\leq (\log
3X)^j.\tag{3.1}$$
In order to estimate the moments in Theorem 2, we have to understand the behavior of the sums
$$\sum_{n=1}^{\infty}
\frac{d_{z_1}(n)d_{z_2}(n)}{n^{2\sigma}}=\prod_{p}\left(\sum_{a=0}^{\infty}
\frac{d_{z_1}(p^a)d_{z_2}(p^a)}{p^{2\sigma a}}\right),$$
for certain $z_1,z_2\in {\Bbb C}$ and $\sigma>1/2$.
The first step is to use the following integral representation for the sum $\sum_{a=0}^{\infty}
d_{z_1}(p^a)d_{z_2}(p^a)/p^{2\sigma a}$.
\proclaim{Lemma 3.2} For all $z_1,z_2\in {\Bbb C}$ and $\sigma>1/2$ we have
$$\sum_{a=0}^{\infty}
\frac{d_{z_1}(p^a)d_{z_2}(p^a)}{p^{2\sigma a}}=\frac{1}{2\pi}\int_{-\pi}^{\pi}\left(1-\frac{e^{i\theta}}{p^{\sigma}}\right)^{-z_1}
\left(1-\frac{e^{-i\theta}}{p^{\sigma}}\right)^{-z_2}d\theta
.$$
\endproclaim
\demo{Proof} We have that
$$
\align
&\frac{1}{2\pi}\int_{-\pi}^{\pi}\left(1-\frac{e^{i\theta}}{p^{\sigma}}\right)^{-z_1}
\left(1-\frac{e^{-i\theta}}{p^{\sigma}}\right)^{-z_2}d\theta=\frac{1}{2\pi}\int_{-\pi}^{\pi}\sum_{a=0}^{\infty}
\frac{d_{z_1}(p^a)e^{i\theta a}}{p^{\sigma a}}\sum_{b=0}^{\infty}\frac{d_{z_2}(p^b)e^{-i\theta b}}{p^{\sigma b}}d\theta\\
&= \sum_{a,b\geq 0}
\frac{d_{z_1}(p^a)d_{z_2}(p^b)}{p^{\sigma(a+b)}}\frac{1}{2\pi}\int_{-\pi}^{\pi}e^{i(a-b)\theta}d\theta=\sum_{a=0}^{\infty}
\frac{d_{z_1}(p^a)d_{z_2}(p^a)}{p^{2\sigma a}}.
\endalign$$
\enddemo
\noindent Using this Lemma we will prove an upper bound for the following sum over the divisor function $d_k(n)$. This will be one of the ingredients to prove Theorem 2.
\proclaim{ Lemma 3.3} Let $k$ be a large real number. For any $0\leq \alpha\leq \frac{\log_3 k}{2\log k}$ we have
$$ \sum_{n=1}^{\infty}\frac{d_k^2(n)}{n^{2-\alpha}}\leq \exp(2k\log\log k(1+o(1))).$$
\endproclaim
\demo{Proof} Let $\sigma=1-\alpha/2$, and put $r=(2k)^{1/\sigma}$. By Lemma 3.2 we have that
$$ \sum_{n=1}^{\infty}\frac{d_k^2(n)}{n^{2-\alpha}}= \prod_p\left(\frac{1}{2\pi}\int_{-\pi}^{\pi}\left(1-\frac{e^{i\theta}}{p^{\sigma}}\right)^{-k}
\left(1-\frac{e^{-i\theta}}{p^{\sigma}}\right)^{-k}d\theta\right).
$$
Now for $p>\sqrt{r}$ (which means that $p^{2\sigma}>2k$) we have that
$$
\align
\frac{1}{2\pi}\int_{-\pi}^{\pi}\left(1-\frac{e^{i\theta}}{p^{\sigma}}\right)^{-k}
\left(1-\frac{e^{-i\theta}}{p^{\sigma}}\right)^{-k}d\theta
&=\frac{1}{2\pi}\int_{-\pi}^{\pi}\exp\left(\frac{2k}{p^{\sigma}}\cos\theta+O\left(\frac{k}{p^{2\sigma}}\right)\right)d\theta\\
&= I_0\left(\frac{2k}{p^{\sigma}}\right)\left(1+O\left(\frac{k}{p^{2\sigma}}\right)\right).\tag{3.2}\\
\endalign
$$
Moreover one can see that $ \frac{1}{2\pi}\int_{-\pi}^{\pi}\left(1-\frac{e^{i\theta}}{p^{\sigma}}\right)^{-k}
\left(1-\frac{e^{-i\theta}}{p^{\sigma}}\right)^{-k}d\theta\leq \left(1-\frac{1}{p^{\sigma}}\right)^{-2k}.$
Hence combining this trivial bound with equation (3.2), and using Lemma 3.1 and Lemma 2.1 we deduce that
$$\align
\sum_{n=1}^{\infty}\frac{d_k^2(n)}{n^{2-\alpha}}&\leq e^{o(k)}\prod_{p\leq \sqrt{r}}\left(1-\frac{1}{p^{\sigma}}\right)^{-2k}\prod_{p>\sqrt{r}}I_0\left(\frac{2k}{p^{\sigma}}\right)\\
&\leq \exp\left(2k\sum_{p\leq r}\frac{1}{p^{\sigma}}+O\left(k+\sum_{p<\sqrt{r}}\log k+\sum_{p>r}\frac{k^2}{p^{2\sigma}}\right)\right)\\
&\leq \exp(2k\log\log k(1+o(1)),
\endalign
$$
which completes the proof.
\enddemo
In order to prove Theorem 1, we have to estimate certain complex moments of $L(1,\chi)$ which, by Theorem 2, are asymptotic to  sums over complex divisor functions. Indeed an essential step to estimate the distribution functions $\Psi(\tau)$ and $\Psi_q(\tau)$ is the following result
\proclaim{Proposition 3.4} For large $s$ we have
$$ \sum_{n=1}^{\infty}
\frac{d_{\frac{s}{2i}}(n)d_{-\frac{s}{2i}}(n)}{n^2}= \exp\left(s\log\log s + C_2 s + C_1\frac{s}{\log s}+O\left(\frac{s}{\log^2 s}\right)\right),$$
where $C_1$ and $C_2$ are defined by equations (1) and (3) respectively.

\endproclaim
Using Lemma 3.2 one can see that
$$\sum_{n=1}^{\infty}
\frac{d_{\frac{s}{2i}}(n)d_{-\frac{s}{2i}}(n)}{n^2}=\prod_{p}\left(\frac{1}{2\pi}\int_{-\pi}^{\pi}\left(1-\frac{e^{i\theta}}{p}\right)^{-\frac{s}{2i}}
\left(1-\frac{e^{-i\theta}}{p}\right)^{\frac{s}{2i}}d\theta\right).$$
The first step in proving Proposition  3.4 is to study the function $$ g_s(\theta):=\left(1-\frac{e^{i\theta}}{p}\right)^{-\frac{s}{2i}}\left(1-\frac{e^{-i\theta}}{p}\right)^{\frac{s}{2i}},$$
as a function of $\theta$ on the interval $[-\pi,\pi]$. We prove the following
\proclaim{Lemma 3.5} For $s>0$, the function $g_s(\theta)$ is a real valued positive continuous function of  $\theta$ on $ [-\pi,\pi]$. Furthermore for $\theta_p:=\arccos\left(\frac{1}{p}\right)$ we have

$$\min_{\theta\in[-\pi,\pi]} g_s(\theta)=g_s(-\theta_p)=\exp\left(-s\arctan\left(\frac{1}{\sqrt{p^2-1}}\right)\right),$$
and
$$\max_{\theta\in[-\pi,\pi]} g_s(\theta)=g_s(\theta_p)=\exp\left(s\arctan\left(\frac{1}{\sqrt{p^2-1}}\right)\right).$$
\endproclaim

\demo{Proof} Let
$
 f_s(\theta):=\log g_s(\theta).$ Then
$$
\align
f_s(\theta)&=-\frac{s}{2i}\log\left(1-\frac{e^{i\theta}}{p}\right)
+\frac{s}{2i}\log\left(1-\frac{e^{-i\theta}}{p}\right)\\
&=\frac{s}{2i}\sum_{n=1}^{\infty}\frac{e^{in\theta}}{p^n n}-\frac{s}{2i}\sum_{n=1}^{\infty}\frac{e^{-in\theta}}{p^n n}=s\sum_{n=1}^{\infty}\frac{\sin(n\theta)}{p^n n},
\endalign
$$
which is the Fourier expansion of $f_s(\theta)$. This implies that $f_s(\theta)$ is an odd continuous real valued  function of $\theta$, from which the first assertion of the Lemma follows.
 By the absolute and uniform convergence of the Fourier expansion of $f_s(\theta)$ we have
$$\frac{f_s'(\theta)}{s}=\sum_{n=1}^{\infty}\frac{\cos(n\theta)}{p^n}= \text{Re}\sum_{n=1}^{\infty}\frac{e^{in\theta}}{p^n}=\text{Re}\left(\frac{e^{i\theta}}{p}\frac{1}{1-e^{i\theta}/p}\right)=\frac{p\cos\theta-1}{|p-e^{i\theta}|^2}.$$
Now the roots of the equation $p\cos\theta=1$ are $\pm \theta_p$, and we can see that
 $f_s'(0)>0$ and $f_s'(-\pi)=f_s'(\pi)<0$. Therefore using that $f_s(\theta)$ is odd, we deduce that
$$\max_{\theta\in[-\pi,\pi]}f_s(\theta)=f_s(\theta_p), \text{ and } \min_{\theta\in[-\pi,\pi]}f_s(\theta)=f_s(-\theta_p)=-f_s(\theta_p).$$
Now
$$ 1-\frac{e^{i\theta_p}}{p}=1-\frac{\cos\theta_p}{p}-i\frac{\sin\theta_p}{p}=1-\frac{1}{p^2}-i\frac{\sqrt{1-1/p^2}}{p}=\left(1-\frac{1}{p^2}\right)\left(1-\frac{i}{\sqrt{p^2-1}}\right).$$
Similarly
$$ 1-\frac{e^{-i\theta_p}}{p}=\left(1-\frac{1}{p^2}\right)\left(1+\frac{i}{\sqrt{p^2-1}}\right).$$
To conclude the proof we note that
$$
\align
f_s(\theta_p)&=\frac{s}{2i}\left(\log\left(1-\frac{e^{-i\theta_p}}{p}\right)-\log\left(1-\frac{e^{i\theta_p}}{p}\right)\right)\\
&=\frac{s}{2i}\log\left(\frac{1+\frac{i}{\sqrt{p^2-1}}}{1-\frac{i}{\sqrt{p^2-1}}}\right)=s
\arctan\left(\frac{1}{\sqrt{p^2-1}}\right).
\endalign
$$

\enddemo

\demo{Proof of Proposition 3.4}
First we define  $$
  h(t):=\left\{\aligned &\log I_0(t) \quad \quad \text{ if } 0\leq t<1,\\  &\log I_0(t)-t \ \text{ if } 1\leq t. \endaligned \right.
$$
For $s>0$ let
$$ E_p(s):=\frac{1}{2\pi}\int_{-\pi}^{\pi}\left(1-\frac{e^{i\theta}}{p}\right)^{-\frac{s}{2i}}\left(1-\frac{e^{-i\theta}}{p}\right)^{\frac{s}{2i}}d\theta.$$
Then by Lemma 3.2 we know that
$$\sum_{n=1}^{\infty}
\frac{d_{\frac{s}{2i}}(n)d_{-\frac{s}{2i}}(n)}{n^2}=\prod_{p}E_p(s).$$

\noindent {\bf Case 1}. $p>\sqrt{s}$. In this case

$$
\align
E_p(s)&= \frac{1}{2\pi}\int_{-\pi}^{\pi}\exp\left(\frac{s}{2i}\frac{e^{i\theta}}{p}-\frac{s}{2i}\frac{e^{-i\theta}}{p}+O\left(\frac{s}{p^2}\right)\right)d\theta.\\
&=\frac{1}{2\pi}\int_{-\pi}^{\pi}\exp\left(\frac{s}{p}\sin\theta+O\left(\frac{s}{p^2}\right)\right)d\theta= I_0\left(\frac{s}{p}\right)\left(1+O\left(\frac{s}{p^2}\right)\right).\\
\endalign
$$

\noindent {\bf Case 2}. $p\leq\sqrt{s}$. First by Lemma 3.5 we have
$$E_p(s)\leq \exp\left(s\arctan\left(\frac{1}{\sqrt{p^2-1}}\right)\right).\tag{3.3}$$
Furthermore let $\epsilon>0$ be a small number to be chosen later. Then using Lemma 3.5 again gives
$$ E_p(s)\geq \frac{1}{2\pi}\int_{\theta_p-\epsilon}^{\theta_p+\epsilon}\left(1-\frac{e^{i\theta}}{p}\right)
^{-\frac{s}{2i}}\left(1-\frac{e^{-i\theta}}{p}\right)^{\frac{s}{2i}}d\theta.$$ Now if $|\theta-\theta_p|\leq \epsilon$ then one can see that
$$ g_s(\theta)
=\left(1-\frac{e^{i\theta_p}+O(\epsilon)}{p}\right)
^{-\frac{s}{2i}}\left(1-\frac{e^{-i\theta_p}+O(\epsilon)}{p}\right)^{\frac{s}{2i}}=g_s(\theta_p)\left(1+O\left(\frac{s\epsilon}{p}\right)\right).$$
Now choose $\epsilon=\frac{1}{Bs}$ where $B$ is suitably large constant. This implies that for $|\theta-\theta_p|\leq \epsilon$ we have $ g_s(\theta)\geq \frac{1}{2}g_s(\theta_p),$ which gives
$$ E_p(s)\gg \frac{1}{s}
\exp\left(s\arctan\left(\frac{1}{\sqrt{p^2-1}}\right)\right).\tag{3.4}$$
Therefore if $p\leq \sqrt{s}$ then  the bounds (3.3) and (3.4) imply that
$$ \log E_p(s)= s\arctan\left(\frac{1}{\sqrt{p^2-1}}\right) +O(\log s).$$
Thus we deduce from the cases 1 and 2 that
$$ \sum_{p}\log E_p(s)=\sum_{p\leq\sqrt{s}}s\arctan\left(\frac{1}{\sqrt{p^2-1}}\right) + \sum_{p>\sqrt{s}}\log I_0\left(\frac{s}{p}\right)+ E_1,
$$
where 
$$ E_1 \ll  \sum_{p\leq \sqrt{s}}\log s+\sum_{\sqrt{s}<p}\frac{s}{p^{2}}\ll \frac{s}{\sqrt{s}\log s}+\frac{\sqrt{s}}{\log s}\log s\ll \sqrt{s}.$$
Now by the prime number theorem we have
$$ \sum_{p\leq x}\frac 1p =\log\log x+ c +O\left(\frac{1}{\log^2 x}\right),\tag{3.5}$$
for some constant $c$.
Hence by the fact that
$ \sum_{p>x}\arctan\left(\frac{1}{\sqrt{p^2-1}}\right)-\frac 1p= O(1/x^2),$ we deduce that
$$ \sum_{p\leq x}\arctan\left(\frac{1}{\sqrt{p^2-1}}\right)=\log\log x+C_2+O\left(\frac{1}{\log^2 x}\right).$$
Therefore
$$
\align
\sum_{p}\log E_p(s)&= s\log\log(\sqrt{s})+ C_2s + s\sum_{\sqrt{s}<p\leq s}\frac 1p + \sum_{p>\sqrt{s}} h\left(\frac{s}{p}\right) +O\left(\frac{s}{\log^2 s}\right)\\
&=s\log\log s+ C_2s + \sum_{p>\sqrt{s}} h\left(\frac{s}{p}\right) +O\left(\frac{s}{\log^2 s}\right),\tag{3.6}
\endalign
$$
where the last estimate follows from (3.5). To complete the proof we only need to evaluate the sum over $h(s/p)$. To this end we use the prime number theorem
 in the form
$$\pi(t)=\int_2^t\frac{du}{\log u}+O\left(te^{-8\sqrt{\log t}}\right).\tag{3.7}$$
First by Lemma 3.1
$$ \sum_{p>s^{3/2}}h\left(\frac{s}{p}\right)\ll \sum_{p>s^{3/2}}\frac{s^2}{p^2}\ll \sqrt{s}.\tag{3.8}$$
Now by (3.7) we have
$$ \sum_{\sqrt{s}<p\leq s^{3/2}}h\left(\frac{s}{p}\right)=
\int_{\sqrt{s}}^{s^{3/2}}h\left(\frac{s}{t}\right)d\pi(t)=
\int_{\sqrt{s}}^{s^{3/2}}h\left(\frac{s}{t}\right)\frac{dt}{\log t}+E_2,\tag{3.9}$$
where
$$
E_2\ll h\left(\sqrt{s}\right)\sqrt{s}e^{-4\sqrt{\log s}}+h\left(s^{-1/2}\right)s^{3/2}e^{-4\sqrt{\log s}}+\int_{\sqrt{s}}^{s^{3/2}}\frac{s}{t^2}\left|h'\left(\frac{s}{t}\right)\right|te^{-8\sqrt{\log t}}dt.$$
Now by Lemma 3.1 we can see that $E_2\ll se^{-4\sqrt{\log s}}.$
To estimate the main term we make the change of variables $T=s/t$. Hence we have
$$
\int_{\sqrt{s}}^{s^{3/2}}h\left(\frac{s}{t}\right)\frac{dt}{\log t}=
s\int_{s^{-1/2}}^{s^{1/2}} \frac{h(T)}{T^2\log(s/T)}dT.\tag{3.10}
$$
In the range $s^{-1/2}\leq t\leq s^{1/2}$, we have
$$ \frac{1}{\log(s/t)}=\frac{1}{\log s}\frac{1}{1-\frac{\log t}{\log s}}=\frac{1}{\log s}+O\left(\frac{\log t}{\log^2 s}\right).$$
Therefore
$$ \int_{s^{-1/2}}^{s^{1/2}} \frac{h(t)}{t^2\log(s/t)}dt=\frac{1}{\log s} \int_{s^{-1/2}}^{s^{1/2}} \frac{h(t)}{t^2}dt+O\left(\frac{1}{\log^2 s}\right),\tag{3.11}$$
by Lemma 3.1 using that
$ \int_{0}^{
\infty}\frac{h(t)\log(t)}{t^2}dt\ll 1.$
Finally by appealing to Lemma 3.1 again we get that
$$ \int_{s^{-1/2}}^{s^{1/2}} \frac{h(t)}{t^2}dt=C_1+O\left(\frac{\log s}{\sqrt{s}}\right).\tag{3.12}$$
Therefore from equations (3.8)-(3.12) we deduce that
$$ \sum_{p>\sqrt{s}} h\left(\frac{s}{p}\right)= C_1\frac{s}{\log s}+O\left(\frac{s}{\log ^2 s}\right),$$
which completes the proof.
\enddemo

\head{4. Complex moments of $L(1,\chi)$}\endhead

To prove Theorem 2, an essential step is to show that $L(1,\chi)^z$ can be approximated by a very short Dirichlet polynomial, if $L(s,\chi)$ has no zeros in a wide region inside the critical strip.
\proclaim{ Proposition 4.1}
Assume that $L(s,\chi)$ has no zeros inside the rectangle $\{s:5/8\leq \text{Re}(s)\leq 1 \text{ and } |\text{Im}(s)|\leq 2\log^3q\}$.
Then for $X=q/\log q$ and any complex number $z$ such that $|z|\leq  R(q)$ we have
$$L(1,\chi)^z=\sum_{n=1}^{\infty}\frac{d_z(n)\chi(n)}{n}e^{-n/X}+O\left(\exp\left(-\frac{\log q\log_4q}{8\log_2 q}\right)\right).$$
\endproclaim

\demo{Proof} Since $\frac{1}{2\pi i}\int_{2-i\infty}^{2+i\infty}y^s\Gamma(s)ds= e^{-1/y}$ then
$$\frac{1}{2\pi i}\int_{2-i\infty}^{2+i\infty}L(1+s,\chi)^z\Gamma(s)X^sds= \sum_{n=1}^{\infty}\frac{d_z(n)\chi(n)}{n}e^{-n/X}.$$
we shift the contour to $\Cal{C}$, where $\Cal{C}$ is the path joining
$$ -i\infty, -i(\log q)^3, -\eta-i(\log q)^3, -\eta+i(\log q)^3, +i(\log q)^3, +i\infty,$$
where $\eta=\log_4 q/(4\log_2 q)$.
We encounter a simple pole at $s=0$ which leaves the residue $L(1,\chi)^z$. Using Lemma 2.2 and Stirling's formula we get
$$
\frac{1}{2\pi i}\left(\int_{-i\infty}^{-i(\log q)^3}+\int_{i(\log q)^3}^{+i\infty}\right)L(1+s,\chi)^z\Gamma(s)X^s ds\ll \int_{(\log q)^3}^{\infty}e^{O(|z|\log_2 qt)}e^{-\frac{\pi}{3}t}dt\ll \frac{1}{q}.
$$
Finally using Corollary 2.5 along with Stirling's formula and the fact that $\Gamma(s)$ has a simple pole at $s=0$, we deduce that
$$
\align
&\frac{1}{2\pi i}\left(\int_{-i(\log q)^3}^{-\eta-i(\log q)^3}+\int_{-\eta-i(\log q)^3}^{-\eta+i(\log q)^3}+\int_{-\eta+i(\log q)^3}^{i(\log q)^3}\right)L(1+s,\chi)^z\Gamma(s)X^s ds\\
&\ll \exp\left(-\frac{\pi}{3}\log^3 q+O(|z|\log_3 q)\right)+\frac{1}{\eta}\log^3q X^{-\eta}\exp\left(|z|\log_3q+O(|z|)\right)\\
&\ll \exp\left(-\frac{\log q\log_4q}{8\log_2 q}\right).
\endalign
$$
\enddemo

\demo{Proof of Theorem 2}
Let $S_q^{+}$ be the set of characters $\chi$ such that $L(s,\chi)$ has no zeros in the rectangle $\{s:5/8\leq \text{Re}(s)\leq 1 \text{ and } |\text{Im}(s)|\leq 2\log^3q\}$, and denote by $S_q^{-}$ the complementary subset $S_q \setminus S_q^{+}$. Then by the zero density result (2.2) we know that $$|S_q^{-}|\ll q^{10/11}.\tag{4.1}$$ Our goal is to evaluate
$$ M_q(z_1,z_2)= \frac{1}{\phi(q)}\sum_{\chi\in S_q}L(1,\chi)^{z_1}L(1,\overline{\chi})^{z_2}.$$
The strategy is as follows: we split the summation into two parts, the first over the characters of $S_q^{+}$ and the second over those in $S_q^{-}$. The latter sum can be trivially bounded using (4.1) and Lemma 2.2. For $\chi\in S_q^{+}$, Proposition 4.1 shows that both $L(1,\chi)^{z_1}$ and $L(1,\overline{\chi})^{z_2}$ can be approximated by very short Dirichlet polynomials. Finally we average the corresponding short sums over all characters (the contribution of the characters $\chi\notin S_q^{+}$ being negligible) and use the orthogonality relations to compute the main term.

We have
$$ M_q(z_1,z_2)=M_q^{+}(z_1,z_2)+M_q^{-}(z_1,z_2),$$
where
$$ M_q^{\pm}(z_1,z_2)= \frac{1}{\phi(q)}\sum_{\chi\in S_q^{\pm}}L(1,\chi)^{z_1}L(1,\overline{\chi})^{z_2}.$$
By (4.1) and Lemma 2.2 we have that
$$ M_q^{-}(z_1,z_2)\ll q^{-1/11}\exp(O(\log_2 q(|z_1|+|z_2|)))\ll q^{-1/12}.$$
 Let $X=q/\log q$ and put $k=\max\{[|z_1|]+1,
[|z_2|]+1\}$. By Proposition 4.1 we have that
$$
M_q^{+}(z_1,z_2)= \sum_{n,m\geq
1}\frac{d_{z_1}(n)d_{z_2}(m)}{nm}e^{-(m+n)/X}\frac{1}{\phi(q)}\sum_{\chi\in S_q^{+}}
\chi(n)\overline{\chi(m)} +E_3,\tag{4.2}
$$
where
$$E_3\ll \exp\left(k\log_3q(1+o(1))-\frac{\log q\log_4q}{8\log_2 q}\right)\ll \exp\left(-\frac{\log q\log_4q}{50\log_2 q}\right),$$
by Corollary 2.5. We now extend the main term of the RHS of (4.2) so as to include all
characters $(\hbox{mod } q)$. To this end we use (3.1) and (4.1) to estimate the contribution of the characters
$\chi\notin S_q^{+}$. Indeed this contribution is bounded by
$$ \frac{|S_q^{-}|+2}{\phi(q)}\left(\sum_{n\geq
1}\frac{d_k(n)}{n}e^{-n/X}\right)^2\ll \frac{(\log 3 X)^{2k}}{q^{1/11}}\ll q^{-1/12}.$$
Therefore using the orthogonality relations of characters we deduce that
$$ M_q(z_1,z_2)= \sum\Sb n,m\geq
1\\ (mn,q)=1\\ m\equiv n \text{ mod } q \endSb\frac{d_{z_1}(n)d_{z_2}(m)}{nm}e^{-(m+n)/X}+ O\left(\exp\left(-\frac{\log q\log_4q}{50\log_2 q}\right)\right).\tag{4.3}$$
First we estimate the contribution of the diagonal terms $m=n$. We know that for all $\alpha>0$ we have $1-e^{-t}\leq 2t^{\alpha}$ for all $t>0$. Then choosing $\alpha= \log_3k/(2\log k)$ the
contribution of these terms is
$$ \sum\Sb n=1\\ (n,q)=1\endSb^{\infty}\frac{d_{z_1}(n)d_{z_2}(n)}{n^2}e^{-2n/X}
= \sum\Sb n=1\\
(n,q)=1\endSb^{\infty}\frac{d_{z_1}(n)d_{z_2}(n)}{n^2} + E_4,\tag{4.4}
$$ where
$$ E_4\ll X^{-\alpha}\sum_{n=1}^{\infty} \frac{d_k(n)^2}{n^{2-\alpha}}\ll X^{-\alpha}\exp(2k\log\log k(1+o(1)))\ll \exp\left(-\frac{\log q\log_4q}{8\log_2 q}\right),$$
 by Lemma 3.3. In the sum on the RHS of (4.4) we remove the constraint $(n,q)=1$ at the cost of an error term bounded by
$$ \ll \frac{k^2}{q^2}\sum_{n=1}^{\infty}\frac{d_k(n)^2}{n^2}\ll \frac{1}{q^2}\exp(2k\log\log k(1+o(1)))\ll \frac{1}{q},$$  using Lemma 3.3, and the fact that $d_k(qn)\leq d_k(q)d_k(n)=kd_k(n).$

Furthermore
the off-diagonal terms $m\neq n$ satisfy $m\equiv n \ (\hbox{mod }
q)$ and $(mn,q)=1$, which imply that $\max (m,n)> q$. Since $X=q/\log q$ we deduce that the
contribution of these terms is bounded by
$$ 2\sum_{n=1}^{\infty}\frac{d_k(n)}{n}e^{-n/X}\left(\sum_{m>q}
\frac{d_k(m)}{m}e^{-m/X}\right)\ll \frac{(\log 3X)^k}{\sqrt{q}}\left(\sum_{m>q}
\frac{d_k(m)}{m}e^{-m/(2X)}\right)\ll q^{-\frac{1}{4}}.$$
Thus we deduce that
$$ M_q(z_1,z_2)= \sum_{n=1}^{\infty}\frac{d_{z_1}(n)d_{z_2}(n)}{n^2}+O\left(\exp\left(-\frac{\log q\log_4q}{50\log_2 q}\right)\right),$$
as desired.
\enddemo

\head{5. Estimating the distribution function} \endhead
\noindent First we shall estimate the Laplace transform of the distribution of $\arg L(1,\chi)$ using our previous estimates for purely imaginary moments of $L(1,\chi)$.
\proclaim{Lemma 5.1} In the range $1\ll s<2R(q)$ we have
$$ \Cal{L}_q(s):=\int_{-\infty}^{\infty}se^{sx}\Psi_q(x)dx=\exp\left(s\log\log s + C_2 s + C_1\frac{s}{\log s}+O\left(\frac{s}{\log^2 s}\right)\right).$$

\endproclaim
\demo{Proof} First since $\arg L(1,\chi)=0$ when $\chi$ is a real character, then 
$$ \Psi_q(\tau)= \frac{1}{\phi(q)}|\{\chi \in S_q : \arg L(1,\chi)>\tau\}|.\tag{5.1}$$
Therefore we have
$$
\align
\La(s)&= \int_{-\infty}^{\infty}se^{sx}\frac{1}{\phi(q)}\sum\Sb\chi \in S_q\\ \arg L(1,\chi)>x\endSb 1 dx= \frac{1}{\phi(q)}\sum_{\chi\in S_q}\int_{-\infty}^{\arg L(1,\chi)}se^{sx}dx\\
&=\frac{1}{\phi(q)}\sum_{\chi\in S_q}e^{s\arg L(1,\chi)}=\frac{1}{\phi(q)}\sum_{\chi\in S_q} L(1,\chi)^{\frac{s}{2i}}L(1,\overline{\chi})^{-\frac{s}{2i}},\tag{5.2}\endalign$$
by the fact that $ \arg L(1,\chi)= \frac{1}{2i}(\log L(1,\chi)-\log L(1,\overline{\chi})).$ Moreover changing the order of the sum and integral in equation (5.2) is justified by the fact that $|\arg L(1,\chi)|\ll \log_2q$ which follows from Corollary 2.3. Finally the result follows by combining Theorem 2 and Proposition 3.4.
\enddemo
\demo{Proof of Theorem 1}
\noindent To estimate $\Psi_q(\tau)$ we use the saddle point method. Let $s$ be the solution to the equation
$$ \tau=\log\log s + C_2 +\frac{C_1+1}{\log s}.\tag{5.3}$$
Let $\epsilon>0$ be a small number to be chosen later and define $$s_1:=s(1+\epsilon), \ s_2:=s(1-\epsilon),\text{ and }
\tau_1:=\tau+\frac{\epsilon}{\log s}, \ \tau_2:=\tau-\frac{\epsilon}{\log s}.$$
To prove Theorem 1 we will first show that for this particular choice of $s$ we have

$$ \La(s)=\int_{\tau_2}^{\tau_1}se^{sx}\Psi_q(x)dx \left(1+O\left(\exp\left(-\frac{s}{\log^2 s}\right)\right)\right).\tag{5.4}$$
To this end we use a variant of Rankin's trick. Indeed since $s-s_2>0$

$$ \int_{-\infty}^{\tau_2}e^{sx}\Psi_q(x)dx\leq \int_{-\infty}^{\tau_2}e^{(s-s_2)(\tau_2-x)+sx}\Psi_q(x)dx\leq e^{\epsilon s\tau_2}\int_{-\infty}^{+\infty}e^{s_2x}\Psi_q(x)dx.$$
Therefore by Lemma 5.1 we have
$$
\align
&\frac{1}{\La(s)}\int_{-\infty}^{\tau_2}se^{sx}\Psi_q(x)dx\leq e^{\epsilon s\tau_2}\frac{s\La(s_2)}{s_2\La(s)}\\
&\leq \exp\left(\epsilon s\tau_2 +s_2\log\log s_2-s\log\log s-\epsilon s C_2 - \epsilon C_1\frac{ s}{\log s} + O\left(\frac{s}{\log^2 s}\right)\right).\tag{5.5}\\
\endalign
 $$
Now $$\epsilon s\tau_2= \epsilon s\log\log s + \epsilon C_2 s+ (\epsilon C_1+\epsilon-\epsilon^2) \frac{s}{\log s},$$ and
$$
\align
s_2\log\log s_2&= s(1-\epsilon) \log(\log s+\log (1-\epsilon))\\
&= (1-\epsilon)s\log\log s + (1-\epsilon)\log(1-\epsilon)\frac{s}{\log s}+O\left(\frac{s}{\log^2 s}\right).\\
\endalign
$$
Hence by inserting these two last estimates in equation (5.5) we deduce that
$$
\align
\frac{1}{\La(s)}\int_{-\infty}^{\tau_2}se^{sx}\Psi_q(x)dx &\leq\exp\left(((1-\epsilon)\log(1-\epsilon)+\epsilon-\epsilon^2)\frac{s}{\log s}+O\left(\frac{s}{\log^2 s}\right)\right)\\
&\leq \exp\left(\left(-\frac{\epsilon^2}{2}+O(\epsilon^3)\right)\frac{s}{\log s}+O\left(\frac{s}{\log^2 s}\right)\right).
\endalign
$$
Now we choose $\epsilon= \frac{A}{\sqrt{\log s}}$, where $A>0$ is a suitably large constant, to get
$$  \frac{1}{\La(s)}\int_{-\infty}^{\tau_2}se^{sx}\Psi_q(x)dx\leq \exp\left(-\frac{s}{\log^2 s}\right).\tag{5.6}$$
Similarly on has
$$ \int_{\tau_1}^{+\infty}e^{sx}\Psi_q(x)dx\leq \int_{\tau_1}^{+\infty}e^{(s_1-s)(x-\tau_1)+sx}\Psi_q(x)dx\leq e^{-\epsilon s\tau_1}\int_{-\infty}^{+\infty}e^{s_1x}\Psi_q(x)dx,$$
and using exactly the same method as before we deduce that
$$  \frac{1}{\La(s)}\int_{\tau_1}^{+\infty}se^{sx}\Psi_q(x)dx\leq \exp\left(-\frac{s}{\log^2 s}\right).\tag{5.7}$$
Therefore by combining inequalities (5.6) and (5.7) we get the estimate (5.4). Now since $\Psi_q(x)$ is a non-increasing function  we have
$$\Psi_q(\tau_1)\int_{\tau_2}^{\tau_1}se^{sx}dx \leq\int_{\tau_2}^{\tau_1}se^{sx}\Psi_q(x)dx \leq \Psi_q(\tau_2)\int_{\tau_2}^{\tau_1}se^{sx}dx.$$
Moreover since
$$ \int_{\tau_2}^{\tau_1}se^{sx}dx= \exp\left(s\tau+O\left(\frac{s}{\log^{3/2}s}\right)\right),$$ then by equation (5.4) and  Lemma 5.1 we deduce that
$$ \Psi_q\left(\tau+\frac{\epsilon}{\log s}\right) \leq \exp\left(-\frac{s}{\log s}+O\left(\frac{s}{\log^{3/2} s}\right)\right)\leq \Psi_q\left(\tau-\frac{\epsilon}{\log s}\right).\tag{5.8}$$
Hence it remains only to solve equation (5.3) in $s$. Indeed we have
$$ e^{\tau-C_2}=\log s\exp\left(\frac{C_1+1}{\log s}\right)= \log s + C_1+1 +O\left(\frac{1}{\log s}\right),$$
and then
$$ s=\exp\left(e^{\tau-C_2}-C_1-1\right)
\left(1+O\left(\frac{1}{e^{\tau}}\right)\right).$$
Thus by inserting these two last estimates in (5.8) we get
$$ \Psi_q(\tau)=\exp\left(-\frac{\exp\left(e^{\tau-C_2}-C_1-1\right)}
{e^{\tau-C_2}}
\left(1+O\left(\frac{1}{e^{\tau/2}}\right)\right)\right).$$
Moreover this last estimate holds
uniformly for $\tau\leq \log_3 q+C_2-o(1),$ using equation (5.3) and Lemma 5.1.

Finally concerning the distribution of the random variable $\arg L(1,X)$, its Laplace transform is given by
$$ 
\align
L_X(s)&=\int_{-\infty}^{\infty}se^{sx}\Psi(x)dx= \ex\left(L(1,X)^{\frac{s}{2i}}\overline{L(1,X)}^{\frac{-s}{2i}}\right)=\sum_{n=1}^{\infty}
\frac{d_{\frac{s}{2i}}(n)d_{-\frac{s}{2i}}(n)}{n^2}\\
&= \exp\left(s\log\log s + C_2 s + C_1\frac{s}{\log s}+O\left(\frac{s}{\log^2 s}\right)\right),
\endalign
$$
by Proposition 3.4. Therefore using exactly the same approach we deduce the same estimate for $\Psi(\tau)$, thus proving the Theorem.
\enddemo

{\bf Acknowledgments}. 
The author is supported by a postdoctoral fellowship from the Natural Sciences and Engineering Research Council of Canada and by the Institute for Advanced Study and the National Science Foundation under agreement No. DMS-0635607.

\enddocument